\documentclass[draft]{amsart}

\usepackage{amssymb}
\usepackage{url}
\usepackage{amscd}

\theoremstyle{plain}

\newtheorem*{maintheorem}{Main Theorem}

\theoremstyle{remark}

\newcommand{\im}{\sqrt{-1}\,}

\newcommand{\bC}{\mathbb{C}}  
\newcommand{\bR}{\mathbb{R}}  
\newcommand{\bZ}{\mathbb{Z}}  
  
\newcommand{\bN}{\mathbb{N}}

\newtheorem{Th}{Theorem}  
\newtheorem{Prop}{Proposition}  
\newtheorem{Cor}{Corollary}  
\newtheorem{Lem}{Lemma}  
\newtheorem{Def}{Definition}  
\newtheorem{Ex}{Example}  
\newtheorem{Rmk}{Remark}

\newcommand{\bt}{\begin{Th}}  
\newcommand{\et}{\end{Th}}  
\newcommand{\bp}{\begin{Prop}}  
\newcommand{\ep}{\end{Prop}}  
\newcommand{\bc}{\begin{Cor}}  
\newcommand{\ec}{\end{Cor}}  
\newcommand{\bl}{\begin{Lem}}  
\newcommand{\el}{\end{Lem}}  
\newcommand{\bd}{\begin{Def}}  
\newcommand{\ed}{\end{Def}}  

\newcommand{\be}{\begin{Ex}}  
\newcommand{\ee}{\end{Ex}}  
\newcommand{\br}{\begin{Rmk}}  
\newcommand{\er}{\end{Rmk}}  
  
\begin{document}

%\date{Received 20xx; received in final form 20xx} 

\title[LcK structures and small deformations of Hopf manifolds]{A note on locally conformally K\"ahler structures and small deformations of Hopf manifolds}

\author{Keizo Hasegawa}

\address{Department of Mathematics, School of Science\\
Osaka University, Toyonaka, Osaka 560-0043, Japan}
\email{hasegawa@math.sci.osaka-u.ac.jp}

\address{Faculty of Education,
Niigata University, Ikarashi-nino-cho, Nishi-ku, 
Niigata 950-2181, Japan}

\email{hasegawa@ed.niigata-u.ac.jp}

%\thanks{Partially supported by grant DGICYT PB98-0618}


\begin{abstract}
A Hopf manifold is a compact complex manifold of which the universal covering is
$\bC^n \backslash \{0\}$. 
In this note we show that any Hopf manifold admits a locally conformally K\"ahler structure (shortly lcK structure), by
constructing a complex analytic family around a Hopf manifold of diagonal type, which admits a lcK potential,
and applying a well known fact (due to Ornea and Verbitsky) that the property of
lcK potential is preserved under a complex analytic family over a sufficiently small parameter space.

\end{abstract}

\maketitle

%%%%%%%%%%%%%%%%%%%%%%%%%%%%%%%%%%%%%%%%%%%%%%%%%%%%%%%%%%%%%%%%%%%%%%%%%
% Macros
%%%%%%%%%%%%%%%%%%%%%%%%%%%%%%%%%%%%%%%%%%%%%%%%%%%%%%%%%%%%%%%%%%%%%%%%%

\newcommand\sfrac[2]{{#1/#2}}

\newcommand\cont{\operatorname{cont}}
\newcommand\diff{\operatorname{diff}}

%%%%%%%%%%%%%%%%%%%%%%%%%%%%%%%%%%%%%%%%%%%%%%%%%%%%%%%%%%%%%%%%%%%%%%%%%

\section{Introduction}
A {\em generalized Hopf manifold} or simply a {\em Hopf manifold} $M$ of complex dimension $n$ is a compact complex manifold
of which the universal covering is $W := \bC^n \backslash \{0\}$, where $n$ is a positive integer ($n \ge 2$). Let $G$ be the covering transformation group
of $M$. Then $G$ acts  properly discontinuously without fixed points as a holomorphic automorphism group on $W$. 
We can see \cite{Hg} that $G$ contains an infinite cyclic subgroup $Z$ of finite index, which defines a contraction map on $W$ (see Sec. 2 for definition).
A {\em primary Hopf manifold} $M$ is a compact complex manifold of which the covering transformation group is an infinite cyclic subgroup $Z$, that is, $M = W/Z$.
It turns out that a Hopf manifold in general has a primary Hopf manifold as a finite normal covering. Note that in the literature 
a primary Hopf manifold is often referred to as a Hopf manifold. The generalized definition above is due to Kodaira \cite{K}.

It is well known that a Hopf manifold is the first example of a non-K\"ahler compact complex manifold, that is, it does not
admit a K\"ahler structure compatible with its complex structure. And this is also the case for a generalized Hopf manifold; that is,
any compact complex manifold whose universal covering is $W$ admit no K\"ahler structures. This is easily seen, for instance, due to the
fact that the diffeomorphism type of a primary Hopf manifold is $S^1 \times S^{2n-1}$. On the other hands,
a Hopf manifold is also known to give an example of  a {\em locally conformally K\"ahler manifold} (or shortly {\em lcK manifold}) (see Sec. 2 for definition),
which is a generalization of K\"ahler manifold. A problem wether any (generalized) Hopf manifold admits a lcK structure
is affirmative, as we see in this note. A (primary) Hopf manifold $W/Z$, where $Z$ is the covering transformation group generated by
a contraction (see Sec. 2 for definition), is called {\em of linear type} if $Z$ is generated by a linear contraction (non-singular linear map on $\bC^n$); and it is called
{\em of diagonal type} if $Z$ is generated by a diagonalizable linear contraction.
It is known that a linear Hopf manifold admit a lcK structure. There is recently a preprint \cite{OV3}
%[L. Ornea1 and M. Verbitsky, {\em Non-linear Hopf manifolds are locally conformally K\"ahler}, arXiv.2202.12398], 
where the authors show that a non-linear Hopf manifold
can be embedded into a linear Hopf manifold, and thus admit a lcK structure. There is a more direct and succinct way to show this result.
In fact, the author (of this note) already made a notice at some workshops
(including the workshop \char96\char96{\em Geometric Structures on Complex Manifolds}", Steklov Institute of  Mathematics, 2011, Moscow) that
there is a complex analytic family over $\bC$ deforming a non-linear Hopf manifold into a linear Hopf manifold, and another complex analytic
family over $\bC$ deforming a linear Hopf manifold to a Hopf manifold of diagonal type such that all fibers over $\bC \backslash \{0\}$  are
biholomorphic. Since a Hopf manifold of diagonal type is of Vaisman type and thus admit a lcK potential (see Sec. 2 for definition), applying the result 
(due to Ornea and Verbitsky) that the property of lcK potential is preserved under a complex analytic family over a sufficiently small parameter space,
we can conclude that any Hopf manifold admits a lcK structure.

In this note, we briefly review basic results relating to lcK structures on Hopf manifolds with concise proofs to them (for the seek of completeness),
leading to a proof for the theorem that any Hopf manifold admits a lcK structure. Note that the construction of lcK structures on Hopf manifolds is
non-trivial and historic as seen in a series of papers \cite{GO, Be, KO, OV2}.

%%%%%%%%%%%%%%%%%%%%%%%%%%%%%%%%%%%%%%%%
\section{Locally conformally K\"ahler strucure}

{\Def A {\em locally conformally K\"ahler structure} (or shortly {\em lcK structure})
on a differentiable manifold $M$ is a Hermitian structure $(h, J)$ on $M$ with its associated fundamental form $\Omega$ satisfying 
$d \,\Omega = \theta \wedge \Omega$ for some closed $1$-form $\theta$ (which is called {\em Lee form}).
A lcK manifold $M$ is of {\em Vaisman type} if its Lee form $\theta$ is parallel w.r.t.
the Levi-Civita connection of $h$; or equivalently, the Lee field $\xi = h^{-1} \theta$ is parallel.}

Note that A lcK structure $\Omega$ is {\em locally conformally K\"ahler}, in the sense that
there is a open covering $\{U_i\}$ of $M$ such that $\Omega_i=e^{-\sigma_i} \,\Omega$
is K\"ahler form on $U_i$ for some functions $\sigma_i$, that is, $d \,\Omega_i = 0$.
The condition $d \,\Omega = \theta \wedge \Omega$ is equivalent to the existence of
a global close $1$-form $\theta$ such that $\theta|U_i = d \sigma_i$.

{\be
A complex projective space $\bC P^1$ is a quotient manifold of $W=\bC^2 -\{O\}$ by
the action of $\bC^*$,
$$\phi_{\lambda}: (z_1,z_2) \rightarrow (\lambda z_1, \lambda z_2) \; (\lambda \in \bC^*).$$
On the other hand, a Hopf surface $S$ is a quotient manifold of $W$ by the
action of $\bZ$
$$\psi^k: (z_1,z_2) \rightarrow (\mu^k z_1, \mu^k z_2) \; (k \in \bZ),$$
for some $\mu \in \bC^*\; (|\mu| > 1)$. 

Since $\Gamma=\{\mu^k \,|\, k \in \bZ\}$ is a discrete subgroup of
$\bC^*$ and $\bC^*/\Gamma$ is a complex torus ${T_{\bC}}^1$,
we see $S$ is a ${T_{\bC}}^1$ bundle over $\bC P^1$.

Consider a $(1,1)$-form on $W$,
$$\omega= -\sqrt{-1} \,(d z_1 \wedge d \overline{z_1} + d z_2 \wedge \overline{z_2}),$$
and put
$$\Omega=\frac{1}{|z_1|^2 + |z_2|^2} \, \omega,$$
then $\Omega$ defines a real $2$-form on $S$. $\Omega$ is not closed, but satisfies
$$d \,\Omega= \theta \wedge \Omega,$$
with
$$\theta= - \frac{1}{|z_1|^2 + |z_2|^2}\, (z_1 d \overline{z_1} + z_2 d \overline{z_2}
+ \overline{z_1} d z_1 +  \overline{z_2} d z_2).$$

For $\psi (X)=\theta (JX)$, we have
$$\psi = \sqrt{-1}\, \frac{1}{|z_1|^2 + |z_2|^2}\, (z_1 d \overline{z_1} + z_2 d \overline{z_2}
- \overline{z_1} d z_1 -  \overline{z_2} d z_2),$$
and if we defines $\overline{\omega}=d \psi$, then
$$\overline{\omega}= \frac{-\sqrt{-1}}{(|z_1|^2 + |z_2|^2)^2}
(|z_2|^2 d z_1 \wedge d \overline{z_1} + |z_1|^2 d z_2 \wedge d \overline{z_2} -
\overline{z_1} z_2 \, d z_1 \wedge d \overline{z_2} - \overline{z_2} z_1 \, d z_2 \wedge d \overline{z_1})$$
which is so called Fubini-Study form. In the affine coordinates $\displaystyle z=\frac{z_2}{z_1}$,
 $\overline{\omega}$ is expressed as
$$\overline{\omega}=\frac{-\sqrt{-1}}{(1 + |z|^2)^2} \, d z \wedge d \overline{z}$$
We can express $\bC P^1$ and Hopf surface $S$ as homogeneous complex manifolds.
Since $G={\rm SL}_2(\bC)$ acts on $\bC P^1$ transitively, we have
$$\bC P^1 = G/B,$$
where $B$ is a Borel subgroup of $G$:
$$B= \left\{ \left(
\begin{array}{@{}cc@{}}
\alpha & \gamma \\
0 & \beta\\
\end{array}
\right) \Bigg| \alpha, \beta \in \bC^*, \alpha \beta=1, \gamma \in \bC \right\}
$$
Consider a subgroup $B_\mu$ of $B$
$$B_\mu = \left\{ \left(
\begin{array}{@{}cc@{}}
\mu^t & \gamma \\
0 & \mu^{-t}\\
\end{array}
\right) \Bigg| \, \mu, \gamma \in \bC, |\mu| > 1,  t \in \bZ \right\}.
$$
Then we have $S=G/B_\mu$, and $B/B_\mu$ is a complex torus $T_\bC^1$.
$S$ is a holomorphic fiber bundle over $\bC P^1$ with fiber $T_\bC^1$.
\ee}

Note that this structure is generalized for a compact homogeneous lcK manifold:
it has a structure of a holomorphic fiber bundle over a complex flag manifold
with fiber $T_\bC^1$ \cite{HK, AHK}.

It is known that lcK structure on $M$ may be defined as a K\"ahler structure $\tilde{\omega}$ on 
the universal covering $\tilde{M}$
on which the the fundamental group $\Gamma$ acts homothetically;
that is, for every $\gamma \in \Gamma$, $\gamma^*\tilde{\omega} = \rho(\gamma) \tilde{\omega}$
holds for some positive constant $\rho(\gamma)$. 

{\Def Let $M$ be a lcK manifold. Suppose that the universal covering $\tilde{M}$ admits a
{\em K\"ahler potential} $\Phi$, which is a real positive function on $\tilde{M}$ such that
$\tilde{\omega}=-\sqrt{-1} \partial \overline{\partial} \Phi$
defines a K\"ahler structure on $\tilde{M}$. If the fundamental group $\Gamma$ acts
homothetically on $\Phi$, then we call $\Phi$ a {\em lcK potential} for $M$.
$\tilde{\omega}$ clearly defines a lcK structure on $M$.
}

{\be
A diagonal Hopf surfaces $S=W /\Gamma_{\mu}$,
where $\Gamma_{\mu}$ is 
generated by a contraction $\psi: (z_1,z_2) \rightarrow
(\mu z_1, \mu z_2)$ with $ |\mu| > 1$ on $W$,
admits a lcK potential
$$\Phi(z_1,z_2)= |z_1|^{2} + |z_2|^{2}.$$
We have a K\"ahler structure
$\tilde{\omega} =  - \sqrt{-1} \,(d\, z_1 \wedge d \, \overline{z_1} +
d\, z_2 \wedge d \, \overline{z_2})$ on $W$ for which
$\tilde{\omega}=-\sqrt{-1} \partial \overline{\partial} \Phi$ holds.
\ee}

%%%%%%%%%%%%%%%%%%%%%%%%%%%%%%%%%%%%%%%%%%%%%%
\section{Structure of Hopf manifolds}

We recall some results \cite{Hg} on a structure theorem of the covering transformation groups of (generalized) Hopf manifolds and
their diffeomorphism types. An analytic automorphism $\phi$ of $W$ may be considered as an analytic automorphism of $\bC^n$ 
which fixes the origin (by Hartogs' theorem). $\phi$ is a {\em contraction} if the sequence $\{\phi^n\}$ converges uniformly to $0$ on
any compact neighborhood of the origin, or equivalently for any relatively compact  neighborhoods $U, V$ of the origin, there exists
$N \in \bN$ such that $\phi^n(\overline{V}) \subset U$ holds for any $n \ge N$.

Let $M=W/G$ be a Hopf manifold, where $G$ is the covering transformation group of $M$ consisting
of analytic automorphisms over $\bC^n$ which fixes the origin $\bold 0$. Then $G$ acts on $W$
properly discontinuously and fixed point free. We can express $G$ as
$$G=H \rtimes Z,$$
where $Z$ is an infinite cyclic group generated by a contraction $\phi$ on $W$, and
$H$ is a finite normal subgroup of $G$ \cite{Hg} . There exists $m \in \bN$ such that for ${Z}'=<\phi'> ,\, \phi'=\phi^m$,
$G'=H \times {Z}'$ (direct product) is a normal subgroup of finite index in $G$; and thus $M=W/G$ has a finite normal covering
$M'=W/G'$, which is diffeomorphic to $S^{2n-1}/H \times S^1$, where $H$ is a finite unitary group acting freely on $S^{2n-1}$,
as shown as a consequence of the deformation arguments we describe in the following.

{\Def A complex analytic family $\pi: {\mathcal M} \rightarrow U$ is a proper holomorphic surbmersion
from a complex manifold ${\mathcal M}$ into a domain $U \subset \bC^k$. Two complex analytic families $\pi: {\mathcal M} \rightarrow U$
and $\pi': {\mathcal M'} \rightarrow U$ over the same domain $U$ are isomorphic if there is an holomorphic
isomorphism $\psi: {\mathcal M} \rightarrow {\mathcal M'}$ commuting with $\pi$ and $\pi'$.
A complex analytic family $\pi: {\mathcal M} \rightarrow U$ is said to be complete at $p \in U$ if any
complex analytic family  $\pi': {\mathcal M'} \rightarrow U' \, (p \in U')$ with $\pi^{-1}(p)=\pi'^{-1}(p)$,
there exists an open neighborhood $V \subset U' \,(p \in V)$ and
a holomorphic map $q: V \rightarrow U$ such that the induced family $\pi': q^{-1} {\mathcal M} \rightarrow V$ by $q$ is isomorphic to
$\pi'|V$; and it is said to be complete if it is complete at any point $p \in U$.
}
 
For an analytic automorphism $\phi$ over $\bC^n$ which fixes the origin $\bold 0$, we denote the linear part of 
$\phi$ (i.e. the Jaconbian matrix $d \phi (\bold 0)$) by $L(\phi)$. It is shown \cite{Hg} that for the covering transformation group $G$ of
a Hopf manifold $W/G$, $L(G)$, as an automorphism group over $W$ is properly discontinuously without fixed point free, defining
a compact complex manifold $W/L(G)$ which is itself a Hopf manifold. This is more or less trivial for an infinite cyclic group $Z$;
and for a general case, it resorts to the Cartan's Uniqueness theorem. 

We can construct a complex analytic family $\pi: {\mathcal M} \rightarrow \bC$ which transforms $M=\pi^{-1}(1)=W/G$ into $\pi^{-1}(0) = W/L(G)$.
In fact, let $T_t, (t \not=0)$ be an analytic automorphism over $W$ defined by 
$$T_t(z_1, z_2,...,z_n) = (t z_1, t z_2,..., t z_n),$$
and set $g_t= T_t g T_t^{-1}$, $G(t)=\{g_t \;| \;g \in G\}$ and $G(0)=L(G)$.
We define for $g \in G$ an analytic automorphism $\tilde{g}$ over $W \times \bC$ as
$$ \tilde{g}: (z, t) \rightarrow  (g_t (z), t),$$
where $z=(z_1, z_2, ..., z_n) \in W$ and $t \in \bC$. Then $\tilde{G}= \{\tilde{g}\, |\, g \in G\}$ is properly discontinuous and fixed-point-free
as an analytic automorphism group over $W \times \bC$ \cite{Hg}. Let ${\mathcal M} = W \times \bC/\tilde{G}$ and ${\pi}: {\mathcal M} \rightarrow \bC$
be the canonical map induced from the $Pr_2: W \times \bC \rightarrow \bC$. Then ${\pi}: {\mathcal M} \rightarrow \bC$ defines a complex analytic
family with $\pi^{-1}(1)= W/G$ and $\pi^{-1}(0)= W/L(G)$.

We can further deform a Hopf manifold $W/L(G)$ to $W/L_0(G)$
with $L_0(G)= L_0(H) \rtimes L_0(Z)$, where $L_0(Z)$
is the contraction generated by a diagonal matrices $d(\alpha_1, \alpha_2,..., \alpha_n)$
with eigenvalues of $\alpha_1, \alpha_2,..., \alpha_n \; (0 <|\alpha_i| < 1, i=1, 2,..., n) $ of $L(\phi)$ and $L_0(H) \subset U(n)$.
In fact, we can assume that $L(\phi)$ is of Jordan form $J(\alpha, n)$.
Let $T_t, (t \not=0)$ be an analytic automorphism over $W$ defined by 
$$T_t(z_1, z_2,...,z_n) = (t ^{n-1} z_1, t^{n-2} z_2,..., t z_{n-1}, z_n),$$
and set $g_t= T_t g T_t^{-1}$, $L(G)(t)=\{g_t \;| \;g \in L(G)\}$. Note that $L(H)(t)$ may be assumed to stay
in $U(n)$ \cite{Hg}. As for the case of the complex analytic family$\{W/G(t)\}$, 
$\{W/LG(t)\}$ defines a complex analytic family  $\pi': {\mathcal M}' \rightarrow \bC$ which transforms
$M'=\pi'^{-1}(1)=W/L(G)$ into $\pi'^{-1}(0) = W/L_0(G)$.

{\be
The complex analytic family constructed above may be considered as a generalization of the following
complex family of Hopf surfaces due to Kodaira \cite{K}. For analytic automorphisms over $W= \bC^2$,
$$g_t: (z_1, z_2) \rightarrow (\alpha z_1 + t z_2, \alpha z_2)\;, 0 < |\alpha| < 1\;,  t \in \bC,$$
we set $Z_t=<g_t>$ and $M_t=W/Z_t$. Then $\{M_t\}$ defines a complex analytic family, where $M_t (t \not=0)$
are all biholomorphic while $M_0$ is not biholomorphic to $M_t (t \not=0)$.
\ee}

{\br
The complete complex analytic families (so-called {\em Kuranishi families}) for Hopf surfaces was
determined by J. Wehler \cite{W}; for Hopf manifolds of diagonal type by C. Borcea \cite{Bo};
and for any Hopf manifolds by A. Haefliger \cite{Hf}. The result that any deformations of Hopf manifolds are
Hopf manifolds (whose universal covering is $W=\bC^n \backslash \{0\}$) is non-trivial and a consequence of the above results.
\er}
 
%\vspace{15pt}

We recall that any Hopf manifold of diagonal type admits a Vaisman structure. 

{\be
A Hopf manifold of diagonal type with eigen values $\lambda_i, i=1,2,..., n$ admits
a locally homogeneous Vaisman structure if all the eigen values $\lambda_i$ have the
same absolute value (which is less than $1$) \cite{HK}.
For simplicity, we consider here the case of Hopf surfaces.
We construct a locally homogeneous
l.c.K. manifolds $\Gamma \backslash {G}/H$ with a discrete
subgroups $\Gamma$ of ${G}=\bR \times U(2)$ and $H=U(1)$.
Let $\Gamma_{p,q}\, (p,q \not=0)$ be a discrete
subgroup of ${G}$:
$$ 
\Gamma_{p,q}=\{(k,\left(
\begin{array}{cc}
e^{\im p k} & 0\\
0 & e^{\im q  k}
\end{array}
\right))
\in \bR \times U(2)\, |\; k \in \bZ
%\in \R \times U(2)\, |\; k \in \Z
\}.
$$
Then $\Gamma_{p,q} \backslash {G}/H$ is biholomorphic to
a Hopf surface $S_{p,q} =W/\Gamma_{\lambda_1,\lambda_2}$, where
$\Gamma_{\lambda_1,\lambda_2}$ is the cyclic group of automorphisms on
$W$ generated by

$$\phi: (z_1,z_2) \longrightarrow (\lambda_1 z_1, \lambda_2 z_2)$$
with $\lambda_1=e^{-r+\im p}, \lambda_2=e^{-r+\im q}, r > 0.$
In fact, we have a homogeneous complex structure $J_r$ on ${G}/H$
for which the biholomorphism $\Phi_r: {G}/H = \bR \times S^3 \rightarrow W$ defined by
$(t, z_1,z_2) \longrightarrow (e^{-r t}  z_1, e^{-r t} z_2)$ induces
a biholomorphism between $\Gamma_{p,q} \backslash {G}/H$ and $S_{p,q}$.
Note that in case $p=q$, $S_{p,q}$ defines a homogeneous Vaisman Hopf surface 
$S^1 \times SU(2)$ \cite{HK}. For the standard homogeneous Sasaki structure $\psi$ (see Example 1)
on $S^3$ and $\theta=d t$ on $\bR$, we have a locally homogeneous Vaisman structure $\Omega$ on $S_{p,q}$:
$$\Omega = - \theta \wedge \psi + d \psi.$$
\ee}

For a Hopf manifold of diagonal type with eigen values $\lambda_i, i=1,2,..., n$, 
where $\lambda_i = e^{-r_i + \sqrt{-1} p_i}$ may have different absolute values $|\lambda_i|= e^{-r_i}, r_i > 0$, 
we consider the biholomorphism $\Phi:  \bR \times S^{2n-1} \rightarrow W$ defined by
$$\Phi: (t, z_1,z_2, ..., z_n) \longrightarrow (e^{-r_1 t}  z_1, e^{-r_2 t} z_2, ..., e^{-r_n t} z_n).$$
For simplicity we consider here again the case of Hopf surfaces $S_{\lambda_1, \lambda_2}$.
In place of the standard Sasaki structure $\psi$,
we take a {\em weighted Sasaki form}:
$$\psi = \sqrt{-1}\, \frac{1}{r_1 |z_1|^2 + r_2 |z_2|^2}\, (z_1 d \overline{z_1} + z_2 d \overline{z_2}
- \overline{z_1} d z_1 -  \overline{z_2} d z_2).$$
The discrete group $\Gamma_{p_1,p_2}\, (p_1,p_2 \not=0)$ acts on $\bR \times S^3$ preserving the weighted Sasaki structure $\psi$,
and $\Phi$ induces a biholomorphism from $\Gamma_{p_1,p_2}  \backslash \bR \times S^{2n-1}$ to $S_{\lambda_1, \lambda_2}$.
For the Lee  form $\theta = dt$ on $\bR$, we have a Vaisman structure $\Omega$ on $S_{\lambda_1, \lambda_2}$:
$$\Omega = - \theta \wedge \psi + d \psi.$$

%%%%%%%%%%%%%%%%%%%%%%%%%%%%%%%%%%%%%%%%%%%%%
\section{Main Theorem}

We now give a proof for the main theorem, based on the above arguments.

%\noindent {\bf Main Theorem}\; 
{\maintheorem 
Any Hopf manifold admits a locally conformally K\"ahler structure.}

%{\it Proof.}
\begin{proof}
We have seen that any Hopf manifold has a primary Hopf manifold as a finite normal covering; and thus
it is sufficient to show the assertion for a primary Hopf manifold. We have just seen that a Hopf manifold of diagonal type
has a Vaisman structure. As is known \cite{OV1}, a compact Vaisman manifold $M$ admits a K\"ahler potential $\Phi$
on the universal covering $\widetilde{M}$ such that the fundamental group $\Gamma$ acts homothetically on $\Phi$.
As observed in \cite{OV1}, for a complex analytic family $\pi: {\mathcal M} \rightarrow U$, where $U$ is an open neighborhood of 
the origin $0$ in $\bC^n$, with $M=\pi^{-1}(0)$, if $M$ admits a a K\"ahler potential $\Phi$, then $M_t= \pi^{-1}(t)$ admit
also a K\"ahler potential for a sufficiently small $t \in U$. This is seen as follows:  $\tilde{\omega}=-\sqrt{-1} \,\partial \overline{\partial} \Phi$
defines a K\"ahler structure on $\tilde{M}$. For $t \in U$, let ${\partial}_t$ and $\overline{\partial}_t$ be partial derivatives with respect to
the complex structures $J_t$ on $M_t$ and also on $\widetilde{M_t}$. For a sufficiently small $t \in U$, 
$\tilde{\omega}=-\sqrt{-1} \,\partial_t \,\overline{\partial}_t \Phi$ is positive definite for a coordinate neighborhood $O_p$ of each point $p \in M$
and thus also for $\widetilde{O_p}$ of $\widetilde{M}$. Since $M$ is compact, we can take a finite number of such evenly covered $O_p$ of $M$.
Hence, $\tilde{\omega}=-\sqrt{-1} \,\partial_t \,\overline{\partial}_t \Phi$ defines a K\"ahler potential for a sufficiently small $t \in U$.
We have seen that any primary Hopf manifold can be deformed through a complex analytic family into a linear primary Hopf manifold;
and any linear primary Hopf manifold can be also deformed into a Hopf manifold of diagonal type. It should be emphasized that the
complex analytic family $\pi: {\mathcal M} \rightarrow M$ we have constructed is such that $\pi^{-1}(t), t \not=0$ are all biholomorphic; and thus
for a given Hopf manifold $M$, $M$ is biholomorphic to $M_t$ for arbitrarily small $t \in U$. We have therefore completed the proof for
the theorem.
\end{proof}

\end{document}